\documentclass[12pt,a4paper,onecolumn]{article}
\usepackage{textcomp}
\usepackage{tikz}
\usepackage{amssymb}
\usepackage{caption}
\usepackage{color}
\usepackage{verbatim}
\usepackage{fancybox}
\usepackage{titlesec} 
\renewcommand\thesection{\arabic{section}} 
\usepackage{setspace}
\usepackage{mathtools}
\usepackage{amsmath}
\usepackage{mathrsfs}
\usetikzlibrary{positioning,fit,calc}
\usepackage[a4paper,left=1in,top=1in,right=1in,bottom=1in,nohead]{geometry} 
\usepackage[hidelinks]{hyperref}
\numberwithin{equation}{section}
\usepackage{array}
\usepackage{bm}
\usepackage{amsthm}
\usepackage[symbol]{footmisc} 
\def\correspondingauthor{\footnote{Corresponding author. Email: williewong088@gmail.com.}}

\tikzset{block/.style={draw,thick,text width=2cm,minimum height=1cm,align=center},
         line/.style={-latex}}
\newcolumntype{P}[1]{>{\centering\arraybackslash}m{#1}} 

\titleformat{\section}[block]{\large\scshape\bfseries}{\thesection.}{1em}{} 
\titleformat{\subsection}[block]{\bfseries}{\thesubsection.}{1em}{} 

\newtheorem{defn}{Definition}[section]
\newtheorem{thm}[defn]{Theorem}

\newtheorem{cor}[defn]{Corollary}
\newtheorem{lem}[defn]{Lemma}

\newtheorem{prob}[defn]{Problem}

\begin{document}
\pagenumbering{arabic}
\begin{center}
    \textbf{\Large On Cross-intersecting Sperner Families}
\vspace{0.1 in} 
    \\{\large W.H.W. Wong\correspondingauthor{}, E.G. Tay}
\vspace{0.1 in} 
\\National Institute of Education\\Nanyang Technological University, Singapore
\end{center}

\begin{abstract}
\noindent Two families $\mathscr{A}$ and $\mathscr{B}$ are said to be cross-intersecting if $A\cap B\neq\emptyset$ for all $A\in\mathscr{A}$ and $B\in\mathscr{B}$. Given two cross-intersecting Sperner families (or antichains) $\mathscr{A}$ and $\mathscr{B}$ on $\mathbb{N}_n$, we prove that $|\mathscr{A}|+|\mathscr{B}|\le 2{{n}\choose{\lceil{n/2}\rceil}}$ if $n$ is odd, and  $|\mathscr{A}|+|\mathscr{B}|\le {{n}\choose{n/2}}+{{n}\choose{(n/2)+1}}$ if $n$ is even. Furthermore, all extremal and almost-extremal families $\mathscr{A}$ and $\mathscr{B}$ are determined.
\end{abstract}
\section{Introduction}
\noindent For $n\in \mathbb{Z}^+$, let $\mathbb{N}_n=\{1,2,\ldots,n\}$ and $2^{\mathbb{N}_n}$ denote its power set. For any integer $k$, $0\le k\le n$, ${{\mathbb{N}_n}\choose{k}}$ denotes the collection of all $k$-sets of $\mathbb{N}_n$. For any family $\mathscr{A}\subseteq 2^{\mathbb{N}_n}$, $\mathscr{A}^{(k)}$ denotes the collection of $k$-sets in $\mathscr{A}$, i.e. $\mathscr{A}^{(k)}=\mathscr{A}\cap{{\mathbb{N}_n}\choose{k}}$. $\mathscr{A}$ is said to be \textit{intersecting} if $X\cap Y\neq \emptyset$ for all $X,Y\in \mathscr{A}$. Two families $\mathscr{A},\mathscr{B}\subseteq 2^{\mathbb{N}_n}$ are said to be \textit{cross $t$-intersecting} if $|A\cap B|\ge t$ for all $A\in\mathscr{A}$ and all $B\in\mathscr{B}$. If $t=1$, we simply say that $\mathscr{A}$ and $\mathscr{B}$ are \textit{cross-intersecting}.
\indent\par A notion closely related to the intersection of sets is the containment of sets. Two subsets $X$ and $Y$ of $\mathbb{N}_n$ are said to be \textit{independent} if $X\not\subseteq Y$ and $Y\not\subseteq X$. If $X$ and $Y$ are independent, we may say that $X$ is independent of $Y$. An \textit{antichain} or \textit{Sperner family} $\mathscr{A}$ on $\mathbb{N}_n$ is a collection of pairwise independent subsets of $\mathbb{N}_n$, i.e., for all $X,Y\in \mathscr{A}$, $X\not\subseteq Y$.
\noindent\par The Erd\H{o}s-Ko-Rado Theorem is central to the study of intersecting family of sets. It was extended by Hilton and Milner to cross-intersecting families of sets.

\begin{thm} (Erd\H{o}s, Ko and Rado \cite{EKR})
Let $n\in\mathbb{Z}^+$ and $\mathscr{A}\subseteq {{\mathbb{N}_n}\choose{k}}$ be an intersecting family for some integer $k\le \frac{n}{2}$. Then, $|\mathscr{A}|\le {{n-1}\choose{k-1}}$.
\end{thm}

\begin{thm} (Hilton and Milner \cite{HAJW MEC})
Let $n\in\mathbb{Z}^+$ and $\mathscr{A}, \mathscr{B}\subseteq {{\mathbb{N}_n}\choose{k}}$ be nonempty cross-intersecting families for some integer $k\le \frac{n}{2}$. Then, $|\mathscr{A}|+|\mathscr{B}|\le {{n}\choose{k}}-{{n-k}\choose{k}}+1$.
\end{thm}
\noindent\par The above results saw extensions in various forms by Frankl \cite{FP}, Frankl and Kupavskii \cite{FP KA 2}, Frankl and Tokushige \cite{FP TN}, and F{\"u}redi \cite{FZ}. For integers $0\le a,b\le n$, they derived upper bounds of $|\mathscr{A}|+|\mathscr{B}|$ for cross-intersecting families $\mathscr{A}\subseteq {{\mathbb{N}_n}\choose{a}}$ and $\mathscr{B}\subseteq {{\mathbb{N}_n}\choose{b}}$. Additionally, Bey \cite{BC}, Matsumoto and Tokushige \cite{MM TN}, Frankl and Kupavskii \cite{FP KA 1}, and Pyber \cite{PL} also obtained upper bounds concerning the product $|\mathscr{A}|\cdot|\mathscr{B}|$. Here, families $\mathscr{A}$ and $\mathscr{B}$, each consisting elements of a fixed size, are also Sperner families too. 
\indent\par Milner obtained what has now become a well-known analogue for an intersecting antichain.
\begin{thm}(Milner \cite{MEC})
Let $n,k\in\mathbb{Z}^+$ and $\mathscr{A}$ be an intersecting Sperner family on $\mathbb{N}_n$. Then, $|\mathscr{A}|\le {{n}\choose{\lfloor \frac{n+k+1}{2}\rfloor}}$.
\end{thm}
\noindent\textbf{Main results and motivation}
\indent\par Parallel to the above results, we shall prove extensions concerning cross-intersecting Sperner families $\mathscr{A}$ and $\mathscr{B}$ on $\mathbb{N}_n$. We prove the maximum possible sum $|\mathscr{A}|+|\mathscr{B}|$ and the uniqueness of extremal families.
\begin{thm}\label{thmA5.1.4}
Let $\mathscr{A}$ and $\mathscr{B}$ be two cross-intersecting antichains on $\mathbb{N}_n$, where $n\in\mathbb{Z}^+$ and $n\ge 3$. Then, 
\begin{align*}
|\mathscr{A}|+|\mathscr{B}|\le {{n}\choose{\lfloor{(n+1)/2}\rfloor}}+{{n}\choose{\lceil{(n+1)/2}\rceil}}
\end{align*}
Furthermore, equality holds if and only if $\{\mathscr{A},\mathscr{B}\}=\{{{\mathbb{N}_n}\choose{\lfloor{(n+1)/2}\rfloor}}, {{\mathbb{N}_n}\choose{\lceil{(n+1)/2}\rceil}}\}$.
\end{thm}
\noindent\par We should mention that Frankl and Wong \cite{FP WHW} and Ou \cite{OY} independently obtained generalisations of Theorem \ref{thmA5.1.4} on cross $t$-intersecting antichains. In this paper, we present a different proof using the classic Sperner operations. It is worthwhile to note that Scott \cite{SAD} also obtained a different proof (a special case where $k=1$) for Milner's theorem using Sperner's operations.
\noindent\par Furthermore, we show that the almost-extremal families are none other than subfamilies of the extremal families.
\begin{thm}\label{thmA5.1.5}
Let $\mathscr{A}$ and $\mathscr{B}$ be two cross-intersecting Sperner families on $\mathbb{N}_n$, where $n\ge 3$ is an odd integer. Then, $|\mathscr{A}|+|\mathscr{B}|=2{{n}\choose{\lceil{n/2}\rceil}}-1$ if and only if $\mathscr{A}={{\mathbb{N}_n}\choose{\lceil{n/2}\rceil}}$, $\mathscr{B}\subset{{\mathbb{N}_n}\choose{\lceil{n/2}\rceil}}$ and $|\mathscr{B}|={{n}\choose{\lceil{n/2}\rceil}}-1$.
\end{thm}
\begin{thm}\label{thmA5.1.6}
Let $\mathscr{A}$ and $\mathscr{B}$ be two cross-intersecting Sperner families on $\mathbb{N}_n$, where $n\ge 4$ is an even integer. Then, $|\mathscr{A}|+|\mathscr{B}|={{n}\choose{n/2}}+{{n}\choose{(n/2)+1}}-1$ if and only if
\\(i) $\mathscr{A}={{\mathbb{N}_n}\choose{n/2}}$, $\mathscr{B}\subset{{\mathbb{N}_n}\choose{(n/2)+1}}$ and $|\mathscr{B}|={{n}\choose{(n/2)+1}}-1$, or
\\(ii) $\mathscr{A}\subset{{\mathbb{N}_n}\choose{n/2}}$, $|\mathscr{A}|={{n}\choose{n/2}}-1$, and $\mathscr{B}={{\mathbb{N}_n}\choose{(n/2)+1}}$.
\end{thm}
\noindent\par The primary motivation of investigating this problem is its intricate connection with optimal orientations of a special family of graphs, known as the $G$ vertex-multiplications. In 2000, Koh and Tay \cite{KKM TEG 8} introduced $G$ vertex-multiplications and extended the results on complete $n$-partite graphs. Koh and Tay \cite{KKM TEG 11} further studied tree vertex-multiplications and Ng and Koh \cite{NKL KKM} investigated cycle vertex-multiplications. Moreover, Wong and Tay \cite{WHW TEG 6A} use the main results in this paper to derive conditions for vertex-multiplications in $\mathscr{C}_0$ and $\mathscr{C}_1$ for all trees of diameter $4$.
\section{Preliminaries}
\noindent Our main tools include the following classical results in extremal set theory.
\begin{thm}(Sperner \cite{SE}) For any $n\in \mathbb{Z}^+$, if $\mathscr{A}$ is a Sperner family on $\mathbb{N}_n$, then $|\mathscr{A}|\le {{n}\choose{\lfloor{n/2}\rfloor}}$. Furthermore, equality holds if and only if all members in $\mathscr{A}$ have the same size, ${\lfloor\frac{n}{2}\rfloor}$ or ${\lceil{\frac{n}{2}}\rceil}$.
\end{thm}
\noindent\par For a family $\mathscr{A}\subseteq {{\mathbb{N}_n}\choose{k}}$, the \textit{shadow} and \textit{shade} of $\mathscr{A}$ are defined as
\begin{align*}
&\Delta \mathscr{A}=\{ X\subseteq \mathbb{N}_n|\ |X|=k-1, X\subset Y\text{ for some } Y\in \mathscr{A}\},\text{ if }k>0,\text{ and }\\
&\nabla \mathscr{A}=\{ X\subseteq \mathbb{N}_n|\ |X|=k+1, Y\subset X\text{ for some }Y\in \mathscr{A}\},\text{ if }k<n 
\end{align*}
respectively.
\indent\par Recall the following elementary inequalities due to Sperner,
\begin{lem} \cite{SE}\label{lemA5.2.2}
Let $\mathscr{A}$ be a collection of $k$-sets of $\mathbb{N}_n$. Then, 
\begin{align}
&|\Delta \mathscr{A}|\ge |\mathscr{A}|\cdot {{n}\choose{k-1}}/{{n}\choose{k}}=|\mathscr{A}|\cdot \frac{k}{n-k+1},\label{eqA4.2.1} \text{ if }k>0, \text{ and }\\
&|\nabla \mathscr{A}|\ge |\mathscr{A}|\cdot {{n}\choose{k+1}}/{{n}\choose{k}}=|\mathscr{A}|\cdot \frac{n-k}{k+1},\label{eqA4.2.2} \text{ if }k<n.
\end{align}
Furthermore, equality holds if and only if $\mathscr{A}=\emptyset$ or $\mathscr{A}={{\mathbb{N}_n}\choose{k}}$.
\end{lem}
\noindent\par Let $\mathscr{F}\subseteq 2^{\mathbb{N}_n}$. Define the top and bottom sizes $t(\mathscr{F})=\max\{|F|\mid F\in\mathscr{F}\}$ and $b(\mathscr{F})=\min\{|F|\mid F\in\mathscr{F}\}$. Following Sperner, let us define two new families obtained from $\mathscr{A}$.
\begin{align*}
\mathscr{F}_\circ=\mathscr{F}-\mathscr{F}^{(t(\mathscr{F}))}\cup\Delta(\mathscr{F}^{(t(\mathscr{F}))}),\\
\mathscr{F}^\circ=\mathscr{F}-\mathscr{F}^{(b(\mathscr{F}))}\cup\nabla(\mathscr{F}^{(b(\mathscr{F}))}).
\end{align*}
\noindent\par The next proposition is well-known.
\begin{lem}\label{lemA4.2.3}
If $\mathscr{F}$ is a nonempty antichain $(\mathscr{F}\neq\{\emptyset\}$ or $\{\mathbb{N}_n\})$, then both $\mathscr{F}_\circ$ and $\mathscr{F}^\circ$ are antichains and $t(\mathscr{F}_\circ)=t(\mathscr{F})-1$, $b(\mathscr{F}^\circ)=b(\mathscr{F})+1$.
\end{lem}
\noindent\par The above bounds for the shadow $\Delta \mathscr{A}$ and shade $\nabla \mathscr{A}$ are not tight except for $\mathscr{A}=\emptyset$ or $\mathscr{A}={{\mathbb{N}_n}\choose{k}}$. A tight lower bound is given by the celebrated Kruskal-Katona Theorem (KKT). KKT is closely related to the \textit{squashed order} of the $k$-sets. The squash relations $\le_s$ and $<_s$ are defined as follows. For $A,B\in {{\mathbb{N}_n}\choose{k}}$, $A\le_s B$ if the largest element of the symmetric difference $A+B=(A-B)\cup (B-A)$ is in $B$. Furthermore, denote $A <_s B$ if $A\le_s B$ and $A\neq B$. For example, the $3$-subsets of $\mathbb{N}_5$ in squashed order are: $\bm{123} <_s \bm{124} <_s \bm{134} <_s \bm{234} <_s \bm{125} <_s \bm{135} <_s \bm{235} <_s \bm{145} <_s \bm{245} <_s \bm{345}$. Here, we omit the braces and write $\bm{abc}$ to represent the set $\{a,b,c\}$, if there is no ambiguity.
\indent\par We shall denote the collections of the first $m$ and last $m$ $k$-subsets of $\mathbb{N}_n$ in squashed order by $F_{n,k}(m)$ and $L_{n,k}(m)$ respectively. KKT says that the shadow of a family $\mathscr{A}$ of $k$-sets is at least the size of the shadow of the first $|\mathscr{A}|$ $k$-sets in squashed order.
\begin{thm}(Kruskal \cite{KJB}, Katona \cite{KGOH}, and Clements and Lindstr{\"o}m \cite{CGF LB})
~\\Let $\mathscr{A}$ be a collection of $k$-sets of $\mathbb{N}_n$ and suppose the $k$-binomial representation of $|\mathscr{A}|$ is
\begin{align}
|\mathscr{A}|={{a_k}\choose{k}}+{{a_{k-1}}\choose{k-1}}+\ldots+{{a_t}\choose{t}},\label{eqA5.2.3}
\end{align}
where $a_k>a_{k-1}>\ldots>a_t\ge t\ge 1$. Then,
\begin{align*}
&|\Delta \mathscr{A}|\ge|\Delta F_{n,k}(|\mathscr{A}|)|={{a_k}\choose{k-1}}+{{a_{k-1}}\choose{k-2}}+\ldots+{{a_t}\choose{t-1}}.
\end{align*}
\end{thm}
\noindent\par Lieby \cite{LP} proved that the shadow of the first $m$ $k$-sets of $\mathbb{N}_n$ in squashed order has the same cardinality as the shade of the last $m$ $(n-k)$-subsets of $\mathbb{N}_n$ in squashed order.
\begin{lem}(Lieby \cite{LP})
For integer $0\le m\le {{n}\choose{k}}$, $|\Delta F_{n,k}(m)|=|\nabla L_{n,n-k}(m)|$.
\end{lem}
\begin{cor} \label{corA5.2.6}
Let $\mathscr{A}$ be a collection of $k$-sets with (\ref{eqA5.2.3}) satisfied and $k=\frac{n}{2}$. Then,
\begin{align*}
|\nabla \mathscr{A}|\ge|\nabla L_{n,\frac{n}{2}}(|\mathscr{A}|)|={{a_k}\choose{k-1}}+{{a_{k-1}}\choose{k-2}}+\ldots+{{a_t}\choose{t-1}}.
\end{align*}
\end{cor}
\section{Proof of Theorem \ref{thmA5.1.4}}
\noindent The proof of Theorem \ref{thmA5.1.4} is almost immediate with Sperner's theorem if $n$ is odd. If $n$ is even, we shall employ Sperner's operations in a manner similar to that used by Sperner himself to prove Sperner's theorem \cite{SE}.
\\
\\\textit{Proof of Theorem \ref{thmA5.1.4}}: Suppose $n$ is odd. For $\mathscr{X}\in\{\mathscr{A},\mathscr{B}\}$, $|\mathscr{X}|\le{{n}\choose{\lfloor{n/2}\rfloor}}$ by Sperner's theorem with equality if and only if $\mathscr{X}={{\mathbb{N}_n}\choose{\lfloor{n/2}\rfloor}}$ or $\mathscr{X}={{\mathbb{N}_n}\choose{\lceil{n/2}\rceil}}$. The cross-intersecting property implies that the equality condition holds if and only if $\mathscr{A}=\mathscr{B}={{\mathbb{N}_n}\choose{\lceil{n/2}\rceil}}$.
\noindent\par Let $n$ be even. Suppose $b(\mathscr{A})<\frac{n}{2}$. Since $|\nabla \mathscr{A}^{(b(\mathscr{A}))}|>|\mathscr{A}^{(b(\mathscr{A}))}|$ by (\ref{eqA4.2.2}) and by Lemma \ref{lemA4.2.3}, we may choose $\mathscr{A}^*\subset \mathscr{A}^\circ$ such that $|\mathscr{A}^*|=|\mathscr{A}|$. Note that $\mathscr{A}^*$ retains the defining properties of $\mathscr{A}$ and $b(\mathscr{A}^*)>b(\mathscr{A})$. By replacing $\mathscr{A}$ with $\mathscr{A}^*$, and repeating the process, we may assume $b(\mathscr{A})\ge \frac{n}{2}$ and similarly, $b(\mathscr{B})\ge \frac{n}{2}$.
\noindent\par Now, suppose $t(\mathscr{A})>\frac{n}{2}+1$. Since $|\Delta \mathscr{A}^{(t(\mathscr{A}))}|>|\mathscr{A}^{(t(\mathscr{A}))}|$ by (\ref{eqA4.2.1}) and by Lemma \ref{lemA4.2.3}, we may choose $\mathscr{A}_*\subset \mathscr{A}_\circ$ such that $|\mathscr{A}_*|=|\mathscr{A}|$. Similarly, $\mathscr{A}_*$ retains the defining properties of $\mathscr{A}$. In particular, each $A\in \mathscr{A}_*$ maintains a nonempty intersection with each $B\in\mathscr{B}$ since $|B|\ge\frac{n}{2}$. Replace $\mathscr{A}$ with $\mathscr{A}_*$ and repeat the process. Perform the same for $\mathscr{B}$. Consequently, we may assume $\mathscr{X}\subseteq {{\mathbb{N}_n}\choose{n/2}}\cup {{\mathbb{N}_n}\choose{(n/2)+1}}$ for $\mathscr{X}\in\{\mathscr{A},\mathscr{B}\}$.
\noindent\par Since there are $\frac{1}{2}{{n}\choose{n/2}}$ sets $\{X,\bar{X}\}$, $X\in{{\mathbb{N}_n}\choose{n/2}}$ and by the cross-intersecting property, we have
\begin{align}
|\mathscr{A}^{(n/2)}|+|\mathscr{B}^{(n/2)}|\le {{n}\choose{n/2}}.\label{eqA5.3.1}
\end{align}
Noting that every element $X$ in $\mathscr{X}^{((n/2)+1)}$ is independent of every element $Y$ in $\mathscr{X}^{(n/2)}$ for $\mathscr{X}\in\{\mathscr{A},\mathscr{B}\}$, it follows from (\ref{eqA4.2.2}) that 
\begin{align}
|\mathscr{X}^{((n/2)+1)}|\le {{n}\choose{(n/2)+1}}-|\nabla \mathscr{X}^{(n/2)}|\le {{n}\choose{(n/2)+1}}-\frac{n}{n+2}|\mathscr{X}^{(n/2)}|. \label{eqA5.3.2}
\end{align}
So,
\begin{align*}
&\ |\mathscr{A}|+|\mathscr{B}| \\
=&\ |\mathscr{A}^{(n/2)}|+|\mathscr{B}^{(n/2)}|+|\mathscr{A}^{((n/2)+1)}|+|\mathscr{B}^{((n/2)+1)}| \\
\le&\ |\mathscr{A}^{(n/2)}|+|\mathscr{B}^{(n/2)}|+{{n}\choose{(n/2)+1}}-\frac{n}{n+2}|\mathscr{A}^{(n/2)}|+{{n}\choose{(n/2)+1}}-\frac{n}{n+2}|\mathscr{B}^{(n/2)}| \\
=&\ (1-\frac{n}{n+2})(|\mathscr{A}^{(n/2)}|+|\mathscr{B}^{(n/2)}|)+2{{n}\choose{(n/2)+1}}\\
\le&\ (1-\frac{n}{n+2}){{n}\choose{n/2}}+2{{n}\choose{(n/2)+1}}\\
=&\ {{n}\choose{n/2}}+{{n}\choose{(n/2)+1}}.
\end{align*}
\noindent\par It is clear that if $\{\mathscr{A},\mathscr{B}\}=\{{{\mathbb{N}_n}\choose{n/2}},{{\mathbb{N}_n}\choose{(n/2)+1}}\}$, then equality holds. Now, assume $|\mathscr{A}|+|\mathscr{B}|={{n}\choose{n/2}}+{{n}\choose{(n/2)+1}}$. From the proof above, it follows that equality must hold in (\ref{eqA5.3.1}) and (\ref{eqA5.3.2}). By the equality condition in (\ref{eqA4.2.2}), $\{\mathscr{A}^{(n/2)},\mathscr{B}^{(n/2)}\}=\{\emptyset, {{\mathbb{N}_n}\choose{n/2}}\}$ which implies $\{\mathscr{A},\mathscr{B}\}=\{{{\mathbb{N}_n}\choose{n/2}},{{\mathbb{N}_n}\choose{(n/2)+1}}\}$.
\qed

\section{Proofs of Theorems \ref{thmA5.1.5}}
\noindent To prove Theorem \ref{thmA5.1.5}, we shall first obtain some useful properties of the binomial difference $D(n,r)$, defined as follows.
\begin{defn}
For any positive integer $r\le n$, define $D(n,r)={{n}\choose{r-1}}-{{n}\choose{r}}$.
\end{defn}

\begin{lem}\label{lemA5.4.2} For any positive integers $r, m$ and $n$, 
\\(a) $D(n,r)\ \substack{>\\=\\<}\ 0 \iff 2r-1\ \substack{>\\=\\<}\ n$.
\\(b) $D(n-1,r-1)+D(n-1,r)=D(n,r)$ for $r \le n-1$.
\\(c) If $n\ge 2r-1$ and $n> m$, then $D(n,r)<D(m,r)$.
\\(d) $D(m,1)\le 0=D(1,1)$ and $D(m,r)\le D(2r-2,r)$ for $r\ge 2$.
\end{lem}
\noindent\textit{Proof}: (a) This follows from $D(n,r)={{n}\choose{r-1}}-{{n}\choose{r}}={{n}\choose{r-1}}(1-\frac{n+1-r}{r})={{n}\choose{r-1}}(\frac{2r-1-n}{r})$.
\\(b) $D(n-1,r-1)+D(n-1,r)=[{{n-1}\choose{r-2}}-{{n-1}\choose{r-1}}]+[{{n-1}\choose{r-1}}-{{n-1}\choose{r}}]=[{{n-1}\choose{r-2}}+{{n-1}\choose{r-1}}]-[{{n-1}\choose{r-1}}+{{n-1}\choose{r}}]
={{n}\choose{r-1}}-{{n}\choose{r}}=D(n,r)$.
\\(c) Note that $D(n+1,r)-D(n,r)=D(n,r-1)<0$ by (b) and (a). So, if $m\ge 2r-1$, then $D(n,r)<D(m,r)$. If $m< 2r-1$, then $D(m,r)>0\ge D(n,r)$ by (a).
\\(d) The first fact is easy to check. If $m\ge 2r-1$, then $D(m,r)\le 0<D(2r-2,r)$ by (a). If $m\le 2r-2$, then by (b) and (a) respectively, $D(m,r)-D(m-1,r)=D(m-1,r-1)\ge 0$.
\qed

\begin{lem}\label{lemA5.4.3} For any positive integers $i$,$j$ and $r$ such that $j\ge2$, $r\le j-1$ and $r\le i\le j-2+r$, $D(i,r)\ge D(j-2+r,r)$.
\end{lem}
\noindent\textit{Proof}: If $r=j-1$, then $i\le j-2+r=2r-1$. So, $D(i,r)\ge 0=D(j-2+r,r)$ by Lemma \ref{lemA5.4.2}(a). If $i=j-2+r$, then we are done. If $r<j-2$ and $i<j-2+r$, then $j-2+r>2r-1$ and thus $D(i,r)>D(j-2+r,r)$ by Lemma \ref{lemA5.4.2}(c).
\qed

\indent\par An identity we will use often in our proofs is Chu Shih-Chieh's identity (CSC), which is also known as ``Hockey Stick Identity''. (See \cite{KKM TEG 0} for more details.)
\begin{lem} (Chu Shih-Chieh's identity)
For any $r,k\in \mathbb{N}$, ${{r}\choose{0}}+{{r+1}\choose{1}}+\ldots+{{r+k}\choose{k}}={{r+k+1}\choose{k}}$.
\end{lem}

\noindent\par The next lemma is an analogue of CSC in terms of $D(n,r)$.
\begin{lem} \label{lemA5.4.5}
For any positive integer $j\ge 2$, $\sum\limits_{r=1}^j{D(j-2+r,r)}=1$.
\end{lem}
\noindent\textit{Proof}: $\sum\limits_{r=1}^j{D(j-2+r,r)}=\sum\limits_{r=0}^{j-1}{{j-1+r}\choose{r}}-\sum\limits_{r=0}^{j-1}{{j-1+r}\choose{r+1}}={{2j-1}\choose{j-1}}-[{{2j-1}\choose{j}}-{{j-2}\choose{0}}]=1$, where we used CSC in the second equality.
\qed

\indent\par Note that Lemmas \ref{lemA5.4.2}, \ref{lemA5.4.3} and \ref{lemA5.4.5} hold for all positive integers while the next two lemmas hold for odd integers $n\ge 3$.
\begin{lem}\label{lemA5.4.6}
Let $n\ge 3$ be an odd integer. Then, $D(i,\lceil\frac{n}{2}\rceil+1)\ge 2$ for $\lceil\frac{n}{2}\rceil+1\le i\le n$.
\end{lem}
\noindent\textit{Proof}: $D(i,\lceil\frac{n}{2}\rceil+1)={{i}\choose{(n+1)/2}}-{{i}\choose{(n+3)/2}}={{i}\choose{(n+1)/2}}[\frac{2((n+3)/2)-1-i}{(n+3)/2}]\ge{{(n+3)/2}\choose{(n+1)/2}}[\frac{2(n+2-i)}{n+3}]= n+2-i\ge 2$.
\qed

\begin{lem}\label{lemA5.4.7}
Let $n\ge 3$ be an odd integer. For any positive integer $m\le {{n}\choose{\lceil n/2\rceil+1}}$, $|\Delta F_{n,\lceil\frac{n}{2}\rceil+1}(m)|\ge m+2$.
\end{lem}
\noindent\textit{Proof}: Let $k=\lceil\frac{n}{2}\rceil+1$ and the $k$-binomial representation of $m$ be $m={{a_k}\choose{k}}+{{a_{k-1}}\choose{k-1}}+\ldots+{{a_t}\choose{t}}$, where $a_k>a_{k-1}>\ldots>a_t\ge t\ge 1$. By KKT, $|\Delta F_{n,\lceil\frac{n}{2}\rceil+1}(m)|={{a_k}\choose{k-1}}+{{a_{k-1}}\choose{k-2}}+\ldots+{{a_t}\choose{t-1}}$.
\indent\par Now, $|\Delta F_{n,\lceil\frac{n}{2}\rceil+1}(m)|-m=\sum\limits_{r=t}^k{D(a_r,r)}$. Since $k\le a_k\le n$, it follows from Lemma \ref{lemA5.4.6}  that 
\begin{align}
D(a_k,k)\ge 2. \label{eqA5.4.1}
\end{align}
If $D(a_r,r)\ge 0$ for all $r=t,t+1,\ldots, k-1$, then $|\Delta F_{n,\lceil\frac{n}{2}\rceil+1}(m)|-m\ge 2$. Now, assume $D(a_r,r)<0$ for some integer $r$, $t\le r\le k-1$. Let $s$ be the smallest integer such that $D(a_r,r)>0$ for all $r=s, s+1, \ldots, k$.
\\
\\Claim: $a_{s}=2s-2$ and $s\ge 2$.
\indent\par Since $D(a_{s},s)>0$, it follows from Lemma \ref{lemA5.4.2}(a) that $a_{s}<2s-1$. Suppose $a_{s}\le 2s-3$. Then, $a_{s-1}\le a_{s}-1\le 2(s-1)-2$. By Lemma \ref{lemA5.4.2}(a), $D(a_{s-1},s-1)>0$, which contradicts the minimality of $s$. Moreover, $D(i,1)\le 0$ for all $i$ by Lemma \ref{lemA5.4.2}(d) and the claim follows.

\indent\par Since $2s-2=a_s> a_{s-1}>a_{s-2}>\ldots>a_t$, we have $a_r\le s-2+r$ for $r=t,t+1,\ldots, s-1$. So,
\begin{align*}
|\Delta F_{n,\lceil\frac{n}{2}\rceil+1}(m)|-m&=\sum\limits_{r=t}^{k}{D(a_r,r)}\\
&\ge D(a_k,k)+\sum\limits_{r=t}^{s}{D(a_r,r)}\\
&= D(a_k,k)+D(2s-2,s)+\sum\limits_{r=t}^{s-1}{D(a_r,r)} \\
&\ge 2+D(2s-2,s)+\sum\limits_{r=t}^{s-1}{D(s-2+r,r)}\\
&= 2+\sum\limits_{r=t}^{s}{D(s-2+r,r)}\\
&\ge 2+\sum\limits_{r=1}^{s}{D(s-2+r,r)}\\
&= 3.
\end{align*}
The first inequality is due to $s<k$ and $D(a_r,r)>0$ for $r=s+1,s+2,\ldots, k-1$ (if there are any of such terms). The second inequality is due to (\ref{eqA5.4.1}) and Lemma \ref{lemA5.4.3}. If $t=1$, the third inequality follows immediately. And, if $t>1$, the third inequality follows from $D(s-2+r,r)<0$ for $r=1,2,\ldots, t-1,$ by Lemma \ref{lemA5.4.2}(a). Invoking Lemma \ref{lemA5.4.5} derives the last equality.
\qed

\indent\par Now, we are ready to prove Theorem \ref{thmA5.1.5}.
\\\textit{Proof of Theorem \ref{thmA5.1.5}}: ($\Leftarrow$) It is straightforward to verify this.
\\
\\($\Rightarrow$) From $|\mathscr{A}|+|\mathscr{B}|=2{{n}\choose{\lceil{n/2}\rceil}}-1$, it follows without loss of generality that $|\mathscr{A}|={{n}\choose{\lceil{n/2}\rceil}}$ and $|\mathscr{B}|={{n}\choose{\lceil{n/2}\rceil}}-1$. By Sperner's theorem, either $\mathscr{A}={{\mathbb{N}_n}\choose{\lfloor{n/2}\rfloor}}$ or $\mathscr{A}={{\mathbb{N}_n}\choose{\lceil{n/2}\rceil}}$. Suppose $\mathscr{A}={{\mathbb{N}_n}\choose{\lfloor{n/2}\rfloor}}$. Then, $A\cap B\neq\emptyset$ for all $A\in \mathscr{A}$ and $B\in\mathscr{B}$ implies $|B|\ge\lceil{\frac{n}{2}}\rceil+1$. Using Sperner's operations as in the proof of Theorem \ref{thmA5.1.4}, we replace all elements of $\mathscr{B}$ of size $>\lceil\frac{n}{2}\rceil+1$ (if any) by an equal number of sets of size $\lceil\frac{n}{2}\rceil+1$ using the respective shadows. It follows that $|\mathscr{B}| \le |{{\mathbb{N}_n}\choose{\lceil{n/2}\rceil+1}}|={{n}\choose{\lceil{n/2}\rceil+1}}$. Then, $|\mathscr{A}|+|\mathscr{B}|\le{{n}\choose{\lfloor{n/2}\rfloor}}+{{n}\choose{\lceil{n/2}\rceil+1}}<2{{n}\choose{\lceil{n/2}\rceil}}-1$, a contradiction.
\indent\par Now, consider $\mathscr{A}={{\mathbb{N}_n}\choose{\lceil{n/2}\rceil}}$. Then, $A\cap B\neq\emptyset$ for all $A\in \mathscr{A}$ and $B\in\mathscr{B}$ implies $|B|\ge\lceil\frac{n}{2}\rceil$. Suppose there exists some $B\in \mathscr{B}$ such that $|B|\ge\lceil\frac{n}{2}\rceil+1$ for a contradiction. Using a similar argument as in the proof of Theorem \ref{thmA5.1.4}, replace all elements of $\mathscr{B}$ of size $>\lceil\frac{n}{2}\rceil+1$ (if any) by an equal number of sets of size $\lceil\frac{n}{2}\rceil+1$ using the respective shadows. So, we may assume $|B|=\lceil\frac{n}{2}\rceil$ or $|B|=\lceil\frac{n}{2}\rceil+1$ for all $B\in\mathscr{B}$. Since $\mathscr{B}$ is a Sperner family, $\Delta \mathscr{B}^{(\lceil\frac{n}{2}\rceil+1)}\cap \mathscr{B}^{(\lceil\frac{n}{2}\rceil)}=\emptyset$ and so
\begin{align*}
{{n}\choose{\lceil{n/2}\rceil}}&\ge |\mathscr{B}^{(\lceil\frac{n}{2}\rceil)}|+|\Delta \mathscr{B}^{(\lceil\frac{n}{2}\rceil+1)}|\\
&\ge |\mathscr{B}^{(\lceil\frac{n}{2}\rceil)}|+|\Delta F_{n,\lceil\frac{n}{2}\rceil+1}(|\mathscr{B}^{(\lceil\frac{n}{2}\rceil+1)}|)|\\
&\ge |\mathscr{B}^{(\lceil\frac{n}{2}\rceil)}|+|\mathscr{B}^{(\lceil\frac{n}{2}\rceil+1)}|+2\\
&= |\mathscr{B}|+2,
\end{align*}
where the second and third inequalities follow from KKT and Lemma \ref{lemA5.4.7}  respectively. It follows that $|\mathscr{A}|+|\mathscr{B}|\le 2{{n}\choose{\lceil{n/2}\rceil}}-2$, a contradiction.
\qed

\section{Proof of Theorem \ref{thmA5.1.6}}
\noindent Now, we consider even integers $n\ge 4$ and proceed in a similar outline as the previous section. Specifically, we probe into the difference between the bounds given in KKT and (\ref{eqA4.2.2}), particularly the term $D^*(n,r,k)$ defined as follows.

\begin{defn}\label{defnA5.5.1} For any positive integers $r,n$ and $k$, define
\begin{align*}
D^*(n,r,k)= \left\{
  \begin{array}{@{}ll@{}}
     {{n}\choose{r-1}}-\frac{k}{k+1}{{n}\choose{r}}, & \text{if}\ r\le n, \\
    0, & \text{otherwise}. \\
  \end{array}\right.
\end{align*}
\end{defn}

\begin{lem}\label{lemA5.5.2}
For all positive integers $i,j,k$ such that $k\ge 2$ and $2j-1\le i \le 2k-1$, $D^*(i,j,k)-D^*(i+1,j,k) \ge \frac{1}{2}$.
\end{lem}
\noindent\textit{Proof}: We shall first prove the following claim by induction.
\\
\\Claim: ${{2j-1}\choose{j-1}} \frac{1}{j+1}\ge \frac{1}{2}$ for all positive integers $j$.
\indent\par If $j=1$, then ${{2j-1}\choose{j-1}} \frac{1}{j+1}={{1}\choose{0}} (\frac{1}{2})=\frac{1}{2}$. Assume ${{2j-1}\choose{j-1}} \frac{1}{j+1}\ge \frac{1}{2}$ for some positive integer $j$. Then, 
\begin{align*}
{{2j+1}\choose{j}} \frac{1}{j+2}=\frac{(2j+1)(2j)(j+1)}{(j+1)(j)(j+2)}\Bigg[{{2j-1}\choose{j-1}}\frac{1}{j+1}\Bigg]=\frac{4j+2}{j+2}\Bigg[{{2j-1}\choose{j-1}}\frac{1}{j+1}\Bigg]>1(\frac{1}{2}),
\end{align*}
where the last inequality follows by induction hypothesis. The claim follows.

\indent\par Now,
\begin{align*}
D^*(i,j,k)-D^*(i+1,j,k)&={{i}\choose{j-1}}-{{i+1}\choose{j-1}}-\frac{k}{k+1}\Bigg[{{i}\choose{j}}-{{i+1}\choose{j}}\Bigg]\\
&=\frac{k}{k+1}{{i}\choose{j-1}}-\frac{j-1}{i+2-j}{{i}\choose{j-1}} \\
&\ge {{i}\choose{j-1}}\Bigg[\frac{j}{j+1}-\frac{j-1}{i+2-j} \Bigg] \\
&\ge {{2j-1}\choose{j-1}} \Bigg[ \frac{j}{j+1}-\frac{j-1}{j+1}\Bigg] \\
&={{2j-1}\choose{j-1}}\cdot \frac{1}{j+1} \\
&\ge \frac{1}{2}\ ,
\end{align*}
where the first and second inequalities follow from the facts that $k\ge j>0$ and $i\ge 2j-1>0$ respectively, and the last inequality follows from the claim.
\qed

\indent\par Note that Lemmas \ref{lemA5.5.3}, \ref{lemA5.5.5} and \ref{lemA5.5.6} are analogous to Lemmas \ref{lemA5.4.3}, \ref{lemA5.4.5} and \ref{lemA5.4.7}.
\begin{lem}\label{lemA5.5.3} For any positive integers $i,k$ and $r$ such that $k\ge 2$, $r\le k-1$, and $r\le i\le k-1+r$, $D^*(i,r,k)\ge D^*(k-1+r,r,k)$ respectively.
\end{lem}
\noindent\textit{Proof}: Note that 
\begin{align}
D^*(2r,r,k)<0 \label{eqA5.5.1}
\end{align}
since $D^*(2r,r,k)={{2r}\choose{r-1}}-\frac{k}{k+1}{{2r}\choose{r}}=\frac{r}{r+1}{{2r}\choose{r}}-\frac{k}{k+1}{{2r}\choose{r}}={{2r}\choose{r}}(\frac{r}{r+1}-\frac{k}{k+1})<0$.
\indent\par If $D^*(i,r,k)<0$, then $i\ge 2r$. Otherwise, $D^*(i,r,k)={{i}\choose{r-1}}-\frac{k}{k+1}{{i}\choose{r}}>{{i}\choose{r-1}}-{{i}\choose{r}}=D(i,r)\ge 0$ by Lemma \ref{lemA5.4.2}(a), a contradiction. Then, by Lemma \ref{lemA5.5.2}, 
\begin{align}
D^*(i,r,k)\ge D^*(i+1,r,k)\ge\ldots\ge D^*(k-1+r,r,k) \label{eqA5.5.2}
\end{align}
\noindent\par If $D^*(i,r,k)\ge 0$, then $D^*(i,r,k)\ge 0>D^*(2r,r,k)\ge D^*(k-1+r,r,k)$ by (\ref{eqA5.5.1}) and (\ref{eqA5.5.2}) respectively.
\qed

\begin{cor}\label{corA5.5.4}
$D^*(k-1+r,r,k)<0$ for all positive integers $k\ge 2$ and $r\le k-1$.
\end{cor}
\noindent\textit{Proof}: As mentioned in Lemma \ref{lemA5.5.3}, $D^*(k-1+r,r,k)\le D^*(2r,r,k)<0$.
\qed

\begin{lem}\label{lemA5.5.5}
For all positive integers $k\ge 2$, $\sum\limits_{r=1}^k{D^*(k-1+r,r,k)}\ge \frac{2}{3}$.
\end{lem}
\noindent\textit{Proof}: $\sum\limits_{r=1}^k{D^*(k-1+r,r,k)}=\sum\limits_{r=0}^{k-1}{{k+r}\choose{r}}-\frac{k}{k+1}\sum\limits_{r=0}^{k-1}{{k+r}\choose{r+1}}={{2k}\choose{k-1}}-\frac{k}{k+1}\big[{{2k}\choose{k}}-{{k-1}\choose{0}}\big] =\frac{k}{k+1} \ge \frac{2}{3}$, where the second equality is due to CSC.
\qed

\indent\par Similar to the odd case, the last few lemmas serve to derive Lemma \ref{lemA5.5.6} which is an improved bound for the size of the shade of a collection of $\frac{n}{2}$-sets compared to Lemma \ref{lemA5.2.2}. As we will see in the proof of Theorem \ref{thmA5.1.6}, this improvement eliminates all possible almost-extremal candidates $\{\mathscr{A},\mathscr{B}\}$ except the ones stated in Theorem \ref{thmA5.1.6}.
\begin{lem}\label{lemA5.5.6}
Let $n\ge 6$ be an even integer. For any positive integer $m<{{n}\choose{n/2}}-1$, $|\nabla L_{n,\frac{n}{2}}(m)|>\frac{n}{n+2}(m)+1$.
\end{lem}
\noindent\textit{Proof}: Let $k=\frac{n}{2}\ge 3$ and the $k$-binomial representation of $m$ be $m={{a_k}\choose{k}}+{{a_{k-1}}\choose{k-1}}+\ldots+{{a_t}\choose{t}}$, where $a_k>a_{k-1}>\ldots>a_t\ge t\ge 1$. By Corollary \ref{corA5.2.6}, $|\nabla L_{n,\frac{n}{2}}(m)|={{a_k}\choose{k-1}}+{{a_{k-1}}\choose{k-2}}+\ldots+{{a_t}\choose{t-1}}$.
\\
\\Case 1. $a_k\le 2k-2$.
\indent\par Now, $|\nabla L_{n,\frac{n}{2}}(m)|-m=\sum\limits_{r=t}^k{D(a_r,r)}$. Since $a_k\le 2k-2$, it follows from Lemma \ref{lemA5.4.2}(a) that $D(a_k,k)>0$, which implies $D(a_k,k)\ge 1$ as it is a difference of two integers. If $D(a_r,r)\ge 0$ for all $r=t,t+1,\ldots, k-1$, then $|\nabla L_{n,\frac{n}{2}}(m)|-\frac{nm}{n+2} \ge |\nabla L_{n,\frac{n}{2}}(m)|-m\ge 1$. 
\noindent\par Now, assume $D(a_r,r)<0$ for some integer $r$, $t\le r\le k-1$. Let $s$ be the smallest integer such that $D(a_r,r)>0$ for all $r=s, s+1, \ldots, k$. As in the proof of Lemma \ref{lemA5.4.7}, it can be shown that $a_{s}=2s-2$ and $s\ge 2$.
\indent\par Similarly, we have $a_r\le s-2+r$ for $r=t,t+1,\ldots, s-1$. So, 
\begin{align*}
|\nabla L_{n,\frac{n}{2}}(m)|-m&=\sum\limits_{r=t}^{k}{D(a_r,r)}\\
&\ge D(2s-2,s)+\sum\limits_{r=t}^{s-1}{D(a_r,r)}\\
&\ge D(2s-2,s)+\sum\limits_{r=t}^{s-1}{D(s-2+r,r)}\\
&\ge\sum\limits_{r=1}^{s}{D(s-2+r,r)}\\
&=1.
\end{align*}
The first inequality is due to $D(a_r,r)>0$ for $r=s+1,s+2,\ldots, k$ (if there are any of such terms). The second inequality is due to Lemma \ref{lemA5.4.3}. If $t=1$, the third inequality follows immediately. And, if $t>1$, the third inequality follows from $D(s-2+r,r)<0$ for $r=1,2,\ldots, t-1,$ by Lemma \ref{lemA5.4.2}(a). Invoking Lemma \ref{lemA5.4.5}  obtains the last equality. Lastly, the strict inequality required in the lemma follows since $\frac{n}{n+2}<1$.
\\
\\Case 2. $a_k= 2k-1$.
\indent\par We shall split into two subcases. First, consider the subcase $t=1$, i.e., there are $k$ terms in the $k$-binomial representation of $m$. Then, $a_1<k$. Otherwise, $m=\sum\limits_{r=1}^{k}{{k-1+r}\choose{r}}={{n}\choose{n/2}}-1$, a contradiction. Since $a_r\le k-1+r$ for $r=1,2,\ldots,k$, it follows that

\begin{align*}
|\nabla L_{n,\frac{n}{2}}(m)|-\frac{n}{n+2}(m)&=\sum\limits_{r=1}^{k}{\Big[{{a_r}\choose{r-1}}-\frac{2k}{2k+2}{{a_r}\choose{r}}\Big]}\\
&=\sum\limits_{r=1}^{k} D^*(a_r,r,k)\\
&\ge\sum\limits_{r=1}^{k} D^*(k-1+r,r,k)+D^*(a_1,1,k)-D^*(k,1,k)\\
&\ge\sum\limits_{r=1}^{k} D^*(k-1+r,r,k)+\frac{1}{2}\\
&\ge\frac{2}{3}+\frac{1}{2}\\
&>1.
\end{align*}
The three inequalities follow from Lemmas \ref{lemA5.5.3}, \ref{lemA5.5.2} and \ref{lemA5.5.5} respectively.
\noindent\par Second, consider the subcase $1<t\le k$. Then,
\begin{align}
|\nabla L_{n,\frac{n}{2}}(m)|-\frac{n}{n+2}(m) &=\sum\limits_{r=t}^{k}{\Big[{{a_r}\choose{r-1}}-\frac{2k}{2k+2}{{a_r}\choose{r}}\Big]}\nonumber\\
&=\sum\limits_{r=t}^{k} D^*(a_r,r,k)\nonumber\\
&\ge\sum\limits_{r=t}^{k}D^*(k-1+r,r,k)\nonumber\\
&\ge\sum\limits_{r=2}^{k}D^*(k-1+r,r,k)\nonumber\\
&=\sum\limits_{r=1}^{k}D^*(k-1+r,r,k)-D^*(k,1,k)\nonumber\\
&\ge \frac{2}{3}-(1-\frac{k^2}{k+1})\nonumber\\
&\ge -\frac{1}{3}+\frac{9}{4}\nonumber\\
&>1\nonumber.
\end{align}
The first inequality is due to $a_r\le k-1+r$ for $r=1,2,\ldots,k-1$, and Lemma \ref{lemA5.5.3} while the second inequality is due to $t\ge 2$ and Corollary \ref{corA5.5.4}. The third inequality follows from Lemma \ref{lemA5.5.5}. Finally, the second last inequality follows since $\frac{k^2}{k+1}$ is increasing for $k\ge 3$.
\qed

\indent\par It is easy to verify that Lemma \ref{lemA5.5.6} does not hold for $n=4$ and $m=3$, and the next lemma will make up for this shortfall. The proof is straightforward and therefore left to the reader.
\begin{lem}\label{lemA5.5.7}
Let $\mathscr{A}$ be an antichain on $\mathbb{N}_4$. If there exists $A\in \mathscr{A}$ such that $|A|=1$ or $|A|=3$, then $|\mathscr{A}|\le 4$. Furthermore, equality holds if and only if $\mathscr{A}$ is one of the families ${{\mathbb{N}_4}\choose{1}}$, ${{\mathbb{N}_4}\choose{3}}$, $\{\bm{1,23, 24, 34}\}$ and $\{\bm{12,13, 14, 234}\}$, up to isomorphism.
\end{lem}

\indent\par Now, we are well-equipped to prove Theorem \ref{thmA5.1.6}.
\\\textit{Proof of Theorem \ref{thmA5.1.6}}: 
\\($\Leftarrow$) It is straightforward to verify this.
\\
\\($\Rightarrow$) Case 1. $n=4$.
\indent\par Without loss of generality, we assume $|\mathscr{A}|>|\mathscr{B}|$. Since $|\mathscr{A}|+|\mathscr{B}|=9$, it follows that $|\mathscr{A}|\ge 5$. Furthermore, $\mathscr{A}$ is a Sperner family implies $|\mathscr{A}|\le {{4}\choose{2}}=6$ by Sperner's theorem.
\\
\\Subcase 1.1. $|\mathscr{A}|=6$.
\indent\par By Lemma \ref{lemA5.5.7}, $\mathscr{A}={{\mathbb{N}_4}\choose{2}}$. Then, $|\mathscr{A}|+|\mathscr{B}|=9$ implies $|\mathscr{B}|=3$. Let $B\in \mathscr{B}$. If $|B|\le 2$, then $A\cap B=\emptyset$ for some $A\in \mathscr{A}$, a contradiction to the cross-intersecting property. If $|B|=4$, then $B=\mathbb{N}_4$ implies $|\mathscr{B}|=1$ as $\mathscr{B}$ is a Sperner family, a contradiction to $|\mathscr{B}|=3$. Hence, (i) follows.
\\
\\Subcase 1.2. $|\mathscr{A}|=5$.
\indent\par By Lemma \ref{lemA5.5.7}, $\mathscr{A}\subset{{\mathbb{N}_4}\choose{2}}$. Then, $|\mathscr{A}|+|\mathscr{B}|=9$ implies $|\mathscr{B}|=4$. By Lemma \ref{lemA5.5.7}, $\mathscr{B}={{\mathbb{N}_4}\choose{1}}$, or ${{\mathbb{N}_4}\choose{3}}$, or $\{\bm{1, 23, 24,34}\}$ or $\{\bm{12,13, 14, 234}\}$, up to isomorphism. If $\mathscr{B}={{\mathbb{N}_4}\choose{1}}$, or $\{\bm{1, 23, 24,34}\}$ or $\{\bm{12,13, 14,234}\}$, up to isomorphism, then there exist some disjoint $A\in \mathscr{A}$ and $B\in \mathscr{B}$, a contradiction. Hence, (ii) follows.
\\
\\Case 2. $n\ge 6$.
\indent\par As shown in the proof of Theorem \ref{thmA5.1.4}, there exists a pair of cross-intersecting antichains $\{\tilde{\mathscr{A}}, \tilde{\mathscr{B}}\}$ such that 
\begin{align}
|\tilde{\mathscr{A}}|+|\tilde{\mathscr{B}}|=|\mathscr{A}|+|\mathscr{B}|={{n}\choose{n/2}}+{{n}\choose{(n/2)+1}}-1 \label{eqA5.5.3}
\end{align}
and $\mathscr{X}\subseteq {{\mathbb{N}_n}\choose{n/2}}\cup {{\mathbb{N}_n}\choose{(n/2)+1}}$ for $\mathscr{X}\in\{\tilde{\mathscr{A}},\tilde{\mathscr{B}}\}$. Similar to (\ref{eqA5.3.1}), we have 
\begin{align}
|\tilde{\mathscr{A}}^{(n/2)}|+|\tilde{\mathscr{B}}^{(n/2)}|\le {{n}\choose{n/2}}\label{eqA5.5.4}
\end{align}
by the cross-intersecting property.
\\
\\Case 2.1. $|\tilde{\mathscr{A}}^{(n/2)}|={{n}\choose{n/2}}$.
\indent\par Then, $|\tilde{\mathscr{B}}^{(n/2)}|=0$ by (\ref{eqA5.5.4}) and $|\tilde{\mathscr{B}}^{((n/2)+1)}|={{n}\choose{(n/2)+1}}-1$ by (\ref{eqA5.5.3}). By Sperner's theorem, $\mathscr{A}=\tilde{\mathscr{A}}^{(n/2)}={{\mathbb{N}_n}\choose{n/2}}$. Now, we need to show $t(\mathscr{B})=\frac{n}{2}+1$. Suppose $t(\mathscr{B})\ge \frac{n}{2}+2$. We reduce to the case where $t(\mathscr{B})=\frac{n}{2}+2$ as follows. If $t(\mathscr{B})>\frac{n}{2}+2$, then we may assume the Sperner operations were performed on $\mathscr{B}$ repeatedly as in Theorem \ref{thmA5.1.4} (i.e., replace $\mathscr{B}$ with some antichain $\mathscr{B}_*\subset \mathscr{B}_\circ$ such that $|\mathscr{B}_*|=|\mathscr{B}|$) until an antichain $\mathscr{B}$ with $t(\mathscr{B})=\frac{n}{2}+2$ was obtained. Now if $|\mathscr{B}^{((n/2)+2)}|>\frac{n-2}{6}$, then by (\ref{eqA4.2.1}) $|\Delta\mathscr{B}^{((n/2)+2)}|-|\mathscr{B}^{((n/2)+2)}|\ge |\mathscr{B}^{((n/2)+2)}|(\frac{6}{n-2})>1$, i.e., $|\Delta\mathscr{B}^{((n/2)+2)}|-|\mathscr{B}^{((n/2)+2)}|\ge 2$. If $|\mathscr{B}^{((n/2)+2)}|<\frac{n-2}{6}<\frac{n}{2}+2$, then by KKT
\begin{align*}
|\Delta\mathscr{B}^{((n/2)+2)}|-|\mathscr{B}^{((n/2)+2)}|&\ge \sum\limits_{i=0}^{|\mathscr{B}^{((n/2)+2)}|-1}\Big[{{(n/2)+2-i}\choose{(n/2)+1-i}}-{{(n/2)+2-i}\choose{(n/2)+2-i}}\Big]\\
&\ge{{(n/2)+2}\choose{(n/2)+1}}-{{(n/2)+2}\choose{(n/2)+2}}\\
&\ge \frac{n}{2}+1\\
&\ge 4.
\end{align*}
In both cases, $\mathscr{A}$ and $\mathscr{B}_\circ$ are cross-intersecting antichains of total size more than ${{n}\choose{n/2}}+{{n}\choose{(n/2)+1}}+1$, a contradiction to Theorem \ref{thmA5.1.4}. Therefore, $t(\mathscr{B})=\frac{n}{2}+1$, in which case, we have (i).
\\
\\Case 2.2. $|\tilde{\mathscr{A}}^{(n/2)}|={{n}\choose{n/2}}-1$.
\indent\par Then, $\nabla \tilde{\mathscr{A}}^{(n/2)}={{\mathbb{N}_n}\choose{(n/2)+1}}$ and $\tilde{\mathscr{A}}^{((n/2)+1)}=\emptyset$. By the cross-intersecting property, $\tilde{\mathscr{B}}^{(n/2)}\subseteq \{\bar{X}_0\}$, where $X_0$ is the only $\frac{n}{2}$-set not in $\tilde{\mathscr{A}}^{(n/2)}$. Suppose $\tilde{\mathscr{B}}^{(n/2)}=\{\bar{X}_0\}$. Since every element $B$ in $\tilde{\mathscr{B}}^{((n/2)+1)}$ is independent of $\bar{X}_0$ in $\tilde{\mathscr{B}}^{(n/2)}$, it follows that $|\tilde{\mathscr{B}}^{((n/2)+1)}|\le |{{\mathbb{N}_n}\choose{(n/2)+1}}-\nabla \tilde{\mathscr{B}}^{(n/2)}|={{n}\choose{(n/2)+1}}-\frac{n}{2}$. It follows that $|\tilde{\mathscr{A}}|+|\tilde{\mathscr{B}}|=|\tilde{\mathscr{A}}^{(n/2)}|+|\tilde{\mathscr{A}}^{((n/2)+1)}|+|\tilde{\mathscr{B}}^{(n/2)}|+|\tilde{\mathscr{B}}^{((n/2)+1)}|\le[{{n}\choose{n/2}}-1]+0+1+[{{n}\choose{(n/2)+1}}-\frac{n}{2}]\le {{n}\choose{n/2}}+{{n}\choose{(n/2)+1}}-3$, a contradiction. Hence, $\tilde{\mathscr{B}}^{(n/2)}=\emptyset$ and $|\tilde{\mathscr{B}}^{((n/2)+1)}|={{n}\choose{(n/2)+1}}$. By (\ref{eqA4.2.1}), $\mathscr{B}=\tilde{\mathscr{B}}={{\mathbb{N}_n}\choose{(n/2)+1}}$. It follows from the cross-intersecting property that $b(\mathscr{A})\ge\frac{n}{2}$, i.e., $\mathscr{A}\subset{{\mathbb{N}_n}\choose{n/2}}$, and we have (ii).
\\
\\Case 2.3. $0<|\tilde{\mathscr{A}}^{(n/2)}|<{{n}\choose{n/2}}-1$.
\indent\par By Corollary \ref{corA5.2.6} and Lemma \ref{lemA5.5.6} respectively, $|\nabla \tilde{\mathscr{A}}^{(n/2)}|\ge |\nabla L_{n,\frac{n}{2}}(|\tilde{\mathscr{A}}^{(n/2)}|)|> \frac{n}{n+2}|\tilde{\mathscr{A}}^{(n/2)}|+1$. Since every element $X$ in $\tilde{\mathscr{A}}^{((n/2)+1)}$ is independent of every element $Y$ in $\tilde{\mathscr{A}}^{(n/2)}$, it follows that $|\tilde{\mathscr{A}}^{((n/2)+1)}|\le {{n}\choose{(n/2)+1}}-|\nabla \tilde{\mathscr{A}}^{(n/2)}|<{{n}\choose{(n/2)+1}}-\frac{n}{n+2}|\tilde{\mathscr{A}}^{(n/2)}|-1$. As in the proof of Theorem \ref{thmA5.1.4}, we can derive that $|\tilde{\mathscr{B}}^{((n/2)+1)}|\le {{n}\choose{(n/2)+1}}-\frac{n}{n+2}|\tilde{\mathscr{B}}^{(n/2)}|$. Using these inequalities and (\ref{eqA5.5.4}), we have
\begin{align*}
&\ |\tilde{\mathscr{A}}|+|\tilde{\mathscr{B}}| \\
=&\ |\tilde{\mathscr{A}}^{(n/2)}|+|\tilde{\mathscr{B}}^{(n/2)}|+|\tilde{\mathscr{A}}^{((n/2)+1)}|+|\tilde{\mathscr{B}}^{((n/2)+1)}| \\
<&\ |\tilde{\mathscr{A}}^{(n/2)}|+|\tilde{\mathscr{B}}^{(n/2)}|+{{n}\choose{(n/2)+1}}-\frac{n}{n+2}|\tilde{\mathscr{A}}^{(n/2)}|-1+{{n}\choose{(n/2)+1}}-\frac{n}{n+2}|\tilde{\mathscr{B}}^{(n/2)}| \\
=&\ (1-\frac{n}{n+2})(|\tilde{\mathscr{A}}^{(n/2)}|+|\tilde{\mathscr{B}}^{(n/2)}|)+2{{n}\choose{(n/2)+1}}-1\\
\le&\ (1-\frac{n}{n+2}){{n}\choose{n/2}}+2{{n}\choose{(n/2)+1}}-1\\
=&\ {{n}\choose{n/2}}+{{n}\choose{(n/2)+1}}+{{n}\choose{(n/2)+1}}-\frac{n}{n+2}{{n}\choose{n/2}}-1 \\
=&\ {{n}\choose{n/2}}+{{n}\choose{(n/2)+1}}-1,
\end{align*}
a contradiction. Hence, this case is not possible.
\\
\\Case 2.4. $|\tilde{\mathscr{A}}^{(n/2)}|=0$.
\indent\par Since $|\tilde{\mathscr{B}}|\le {{n}\choose{n/2}}$ by Sperner's theorem, we have either $|\tilde{\mathscr{A}}^{((n/2)+1)}|={{n}\choose{(n/2)+1}}-1$ or $|\tilde{\mathscr{A}}^{((n/2)+1)}|={{n}\choose{(n/2)+1}}$. The former case is similar to Case 2.1 with the roles of $\tilde{\mathscr{A}}$ and $\tilde{\mathscr{B}}$ swapped. In the latter case, $\mathscr{A}=\tilde{\mathscr{A}}={{\mathbb{N}_n}\choose{(n/2)+1}}$ by (\ref{eqA4.2.1}) and $|\tilde{\mathscr{B}}|={{n}\choose{n/2}}-1$ by (\ref{eqA5.5.3}). Note that $\tilde{\mathscr{B}}^{(n/2)}\neq \emptyset$ since $|\tilde{\mathscr{B}}^{((n/2)+1)}|\le {{n}\choose{(n/2)+1}}<{{n}\choose{n/2}}-1$. Now the cases $|\tilde{\mathscr{B}}^{(n/2)}|={{n}\choose{n/2}}-1$ and $0<|\tilde{\mathscr{B}}^{(n/2)}|<{{n}\choose{n/2}}-1$ follow from Cases 2.2 and 2.3 respectively, with the roles of $\tilde{\mathscr{A}}$ and $\tilde{\mathscr{B}}$ swapped.
\qed

\indent\par We end by proposing a direction for further generalisation of the main results. It is known that in Lubell's \cite{LD} proof of Sperner's theorem, a stronger result known as the Lubell-Yamamoto-Meshalkin (LYM) inequality was derived. In this line of thought, we have the following problem.
\begin{prob} Strengthen Theorem \ref{thmA5.1.4} into its LYM form.
\end{prob}
\section*{Acknowledgement}
\noindent The first author would like to thank the National Institute of Education, Nanyang Technological University of Singapore, for the generous support of the Nanyang Technological University Research Scholarship.

\end{document}